\documentclass[english,11pt]{smfart}

\usepackage{latexsym}
\usepackage{epsfig}

\def\cal{\mathcal}
\font\tenbb=msbm10
\def\cC{\hbox{\tenbb C}}
\def\rR{\hbox{\tenbb R}}

\def\tT{\hbox{\tenbb T}}

\newcommand{\Z}{\mathbb{Z}}
\newcommand{\R}{\mathbb{R}}

\newcommand{\T}{\mathbb{T}}

\def\di{\displaystyle}

\newtheorem{defi}{Definition}
\newtheorem{thm}{Theorem}
\newtheorem{lem}{Lemma}

\begin{document}

\setcounter{tocdepth}{2}
\baselineskip 6mm
\title[The Transversality-torsion phenomenon]{Hyperbolicity versus partial-hyperbolicity and the transversality-torsion phenomenon}
\author{Jacky CRESSON}
\address{Equipe de Math\'ematiques de Besan\c{c}on, CNRS-UMR 6623, 16 route de Gray, 25030
Besan\c{c}on cedex, France.}
\email{cresson@math.univ-fcomte.fr}

\author{Christophe GUILLET}
\address{Institut Universitaire de Technologie, 1,All\'ee des Granges Forestier, 71100
Chalon sur Sa{\^o}ne, France.}
\email{Christophe.Guillet@u-bourgogne.fr}

\maketitle

\begin{abstract}
In this paper, we describe a process to create hyperbolicity in
the neighbourhood of a homoclinic orbit to a partially hyperbolic
torus for three degrees of freedom Hamiltonian systems: the transversality-torsion phenomenon.
\end{abstract}

\noindent {\bf keywords}: Hyperbolicity, Partially hyperbolic tori, Hamiltonian systems.

\section{Introduction}

The aim of this paper is to describe a process of creation of hyperbolicity in a {\it partially hyperbolic} context called the
{\it transversality-torsion phenomenon} introduced in (\cite{cr2},\cite{cr3}). This
process comes from the study of instability (Arnold diffusion) for (at least) three degrees of freedom near-integrable Hamiltonian systems
\cite{ar} and more precisely from the derivation of a Smale-Birkhoff theorem (\cite{ea},\cite{cr3}) for transversal homoclinic partially hyperbolic
tori which come along multiple resonances \cite{tr}. Our starting point is
the following conjecture of R.W. Easton (\cite{ea},p.252) about {\it symbolic dynamics} for transversal homoclinic partially hyperbolic
tori: In \cite{ea}, Easton prove the existence of symbolic dynamics in a neighbourhood of a partially hyperbolic torus whose stable and
unstable manifolds intersect transversally (in a given energy manifold containing the torus). This result is obtained under
hypothesis, the most stringent one being on the linear part of the {\it homoclinic map} (see $\S$.\ref{mise} for a definition and
\cite{ea},p.244), called the {\it homoclinic matrix}. However, Easton conjectures (\cite{ea},p.252) that this assumption can
be weakened, or perhaps cancelled. Moreover, the role of all the parameters of the problem (transversality, torsion of the flow
around the torus) is not clear.\\

In (\cite{cr2},\cite{cr3}), we weaken the homoclinic matrix condition, but mainly, we put in evidence a dynamical and
geometrical phenomenon at the origin of the {\it hyperbolic nature} of symbolic dynamics, the {\it transversality-torsion phenomenon}: the transversality of the
stable and unstable manifold of the torus coupled with the torsion of the flow around the torus give rise to a hyperbolic
dynamics in the neighbourhood of the homoclinic connection. Since this first study, the transversality-torsion phenomenon
has been identified and extended by others. We refer in particular to the papers of M. Gidea and C. Robinson (\cite{gr} p.64)
and M. Gidea and R. De La LLave \cite{gl}, dealing with topological methods in dynamics.\\

In this paper, we prove that the transversality-torsion phenomenon observed in a particular case
in \cite{cr3} arises in a generic situation for three degrees of freedom Hamiltonian systems.\\

The plan of the paper is the following: In $\S$.\ref{object}, we define transversal homoclinic partially hyperbolic tori.
In $\S$.\ref{mise}, we state precisely the {\it hyperbolicity problem}, which can be resume as finding the minimal conditions
(about the dynamics on the torus and the geometry of the intersection of the stable and unstable manifolds) in order to
have a {\it homoclinic transition map}\footnote{See $\S$.\ref{setup} for a definition.} hyperbolic. In $\S$.\ref{transtorpheno}, we solve the hyperbolicity problem for three
degrees of freedom Hamiltonian systems, putting in evidence the transversality-torsion phenomenon, i.e. the fundamental role
of the torsion of the flow on the torus and the transversality of the stable and unstable manifolds to induce hyperbolicity of
the transition map.

\section{Transversal homoclinic partially hyperbolic tori}
\label{object}

In this section, we define partially hyperbolic tori following the paper of S. Bolotin and D. Treschev \cite{bt}.

\subsection{Partially hyperbolic tori}

Let $\cal M$ be a $2m$ dimensional symplectic manifold, and $H$ an analytic Hamiltonian defined on $\cal M$.

\begin{defi}
A weakly reducible, diophantine partially hyperbolic torus for $H$ is a torus for which there exists an analytic
symplectic coordinates system, such that the Hamiltonian takes the form
\begin{equation}
\label{normal}
H(\theta ,I,s,u)=\omega .I + \di {1\over 2}\, AI.I +s.M(\theta )u +O_3 (I,s,u) ,
\end{equation}
where $(\theta ,I,s,u)\in \tT^k \times \rR^k \times \rR^{m-k} \times \rR^{m-k}$,
with the symplectic structure $\nu = dI \wedge d\theta +ds\wedge du$, $A$ is a $k\times k$
symmetric constant matrix, $M$ is a definite positive matrix and for all $k\in \Z^n \setminus \{ 0\}$, we have
$$\mid \omega .k\mid \geq \alpha \mid k\mid^{-\beta},\ \alpha ,\beta >0 .$$
If $M$ is a constant matrix, then the partially hyperbolic torus is say to be reducible.
\end{defi}

In (\cite{bt}, Theorem 1, p. 406), S. Bolotin and D. Treschev prove that this "KAM" definition is equivalent to the dynamical one (see \cite{bt},
Definition 1 and 3, p. 402).
Moreover, for $k=1$ and $k=m-1$, the torus is always reducible.\\

In \cite{bt}, Bolotin and Treschev introduce the notion of {\it nondegenerate} hyperbolic torus, which is a condition of
dynamical nature (see \cite{bt}, definition 3, p.402). In the setting of weakly reducible hyperbolic tori, we can use
the following definition which is equivalent to the dynamical one (see \cite{bt}, Proposition 2,p. 404):

\begin{defi}
A weakly reducible hyperbolic torus is nondegenerate if $\mbox{\rm det} A \not= 0$.
\end{defi}

H. Eliasson \cite{el} and L. Niederman \cite{ni}, have proved the following normal form theorem for $m-1$ dimensional tori:

\begin{thm}
\label{elia}
Let $T$ be a $m-1$ dimensional reducible and non-degenerate diophantine partially hyperbolic torus.
There exists an analytic coordinates system $(x,y,z^+,z^- )$ defined in a neighbourhood $V$ of $T$, such that
\begin{equation}
\label{noreli}
H=\omega . y +\lambda z^- z^+ +O_2 (y,z^+ z^- ) .
\end{equation}
\end{thm}

The geometry of the torus can then be easily described (\cite{bt}): it admits analytic stable (resp. unstable) manifold,
denoted by $W^+ (T)$ (resp. $W^- (T)$), and locally defined in $V$ by:
\begin{equation}
\left .
\begin{array}{lll}
W^+ (T) & = & \{ (x,y,z^+ ,z^- ) \in V ,\ y=0, z^- =0 \} , \\
W^- (T) & = & \{ (x,y,z^+ ,z^- ) \in V ,\ y=0, z^+ =0 \} .
\end{array}
\right .
\end{equation}

\subsection{Transversal homoclinic connection}

In the following, we denote by $\cal H$ the energy submanifold of $\cal M$ containing the torus under consideration.
For convenience, a weakly reducible diophantine partially hyperbolic torus will be called a partially hyperbolic
torus.

\begin{defi}
Let $T$ be a $m-1$ dimensional partially hyperbolic torus. We say that $T$
possesses a transversal homoclinic connection if its stable and unstable manifolds intersect transversally in $\cal H$.
\end{defi}

In this paper, we explore the existence of a {\it hyperbolic} dynamics in a neighbourhood of a transversal homoclinic
connection to a partially hyperbolic torus.

\section{The hyperbolicity problem}
\label{mise}

\subsection{Set-up}
\label{setup}

Let $H$ be a $m$ degree of freedom Hamiltonian system. Let $T$ be a $m-1$ dimensional partially hyperbolic torus of $H$
possessing a transversal homoclinic connection along a homoclinic (at least one)
orbit denoted by $\gamma$. We introduce the following notations and terminology:\\

Let $V$ be the Eliasson's normal form domain (\ref{noreli}). There exists (\cite{ma}),
a Poincar\'e section $\cal S$ of $T$ in $V$, and an analytic coordinates systems in $\cal S$,
denoted by $(\phi ,\rho ,s,u)\in \tT^{m-2} \times \rR \times \rR^{m-2} \times \rR$, such that the Poincar\'e map takes
the form
\begin{equation}
\label{poin}
f(\phi ,s,\rho ,u) =(\phi +\omega +\nu \rho ,\lambda s, \rho ,\lambda^{-1} u ) +O_2 (\rho ,s,u) ,
\end{equation}
where $\omega \in \rR^{m-2}$, $\nu \in \rR^{m-2}$, $0<\lambda <1$, $\nu \rho =(\nu_1 \rho_1 ,\dots ,
\nu_{m-2} \rho_{m-2} )$.

We denote by $f_l (\phi ,s,\rho ,u) =(\phi +\omega +\nu \rho ,\lambda s, \rho ,\lambda^{-1} u )$ the linear part of $f$.\\

We say that the torus $T$ is with {\it torsion} if $\nu_i \not= 0$, for $i=1 ,\dots ,m-2$, and {\it without torsion}
otherwise. We note that a torus is with torsion if and only if it is nondegenerate.\\

Let $p^- = (\phi^- ,0,0,u^- ) \in {\cal S}$
and $p^+ =(\phi^+ ,s^+ ,0,0) \in {\cal S}$, be the last (resp.
the first) point of intersection between $\gamma$ and $\cal S$ along the unstable manifold (resp. the stable manifold).
There exists neighbourhoods $V^+$ and $V^-$ in $\cal S$ of $p^+$ and $p^-$ respectively, and a map
$\Gamma :V^- \rightarrow V^+$, called the {\it homoclinic map}, such that $\Gamma (p^- )=p^+$. The homoclinic map is
of the form
$$\Gamma (p^- +z)=p^+ +\Pi . z +O_2 (z),$$
where $\Pi$ is a matrix, called the {\it homoclinic matrix}. We denote by
$$\Gamma_l (p^- +z)=p^+ +\Pi .z.$$

We denote by $D_n =\{ z\in V^+ \mid f_l^n (z)\in V^- \}$ and $D=\bigcup_{n\geq 1} D_n$. We denote by
$\psi : D\rightarrow V^-$, the {\it transverse map} introduced by J\"urgen Moser \cite{mo} and defined by
$$\psi (z)=f^n (z)\ \mbox{\rm if}\ z\in D_n .$$
We denote by $\psi_l (z)=f_l^n (z)\ \mbox{\rm if}\ z\in D_n$.\\

The differential of $f_l$, denoted by $Df_l$ is the matrix
\begin{equation}
Df_l =
\left (
\begin{array}{cccc}
\mbox{\rm Id} & 0 & {\cal V} & 0 \\
0 & \lambda & 0 & 0 \\
0 & 0 & \mbox{\rm Id} & 0 \\
0 & 0 & 0 & \di\lambda^{-1}
\end{array}
\right )
,
\end{equation}
where $\mbox{\rm Id}$ is the $(m-2)\times (m-2)$ identity matrix and $\cal V$ the diagonal matrix with components
$\nu_i$, $i=1,\dots ,m-2$.\\

In the following, we always work in the Poincar\'e section $\cal S$.\\

Let $C=\{ (u,v)\in \rR^{m-1} \times \rR^{m-1} \mid \parallel u\parallel_1 \leq 1 ,\parallel v\parallel_1 \leq 1 \}$. We denote by
$W_{\mu}  :C \rightarrow V^+$, what we call an {\it Easton's window} (or simply window in the following) defined by
$$W_{\mu} (z)=\mu z +p^+.$$
We consider the map $\Delta :C\rightarrow C$, defined by
$$\Delta =(W_{\mu} )^{-1} \circ \Gamma \circ \psi \circ W_{\mu}.$$
We denote by
$$\Delta_l =(W_{\mu} )^{-1} \circ \Gamma_l \circ \psi_l \circ W_{\mu}.$$
We have $\Delta$ and $\Delta_l$ as close as we want in $C^1$-topology when $\mu \rightarrow 0$ (\cite{ea},p.250).\\

The map $\Gamma \circ \psi$ is called the {\it homoclinic transition map}.\\

In the following, we call {\it linear model} a Hamiltonian system possessing a transversal homoclinic partially
hyperbolic torus $T$ such that the preceding maps are linear in a given coordinates systems. Using the Easton's window
$W_{\mu}$, the general case reduces to prove results on the linear model which are stable under small $C^1$ perturbations.
This is the case for example when dealing with symbolic dynamics.

\subsection{The hyperbolicity problem}

We keep the notations and terminology of the previous section. For all matrix $M$, we denote by $\mbox{spec} (M)$
its spectrum. The {\it hyperbolicity problem} can be formulated as follow:\\

{\bf Hyperbolicity Problem} -- {\it Let $H$ be a $m$ degrees of freedom Hamiltonian system. Let $T$ be a $m-1$ dimensional
partially hyperbolic of $H$ possessing a transversal homoclinic connection. Under which conditions on $n$, $\nu$ and $\Pi$ do we have
$$\mbox{\rm spec}(\Pi.Df_l^n )
\cap S^1 =\emptyset,$$
where $S^1 =\{ z\in \cC,\ \mid z\mid =1\}$ is the unit circle in $\cC$.}\\

If $\mbox{\rm spec}(\Pi.Df_l^n ) \cap S^1 =\emptyset$, then for $\mu$ sufficiently small, i.e. in a given neighbourhood
of the homoclinic orbit, the map $\Delta$ is hyperbolic.\\

This problem is difficult as there exists no results about localization of eigenvalues for the product of two
matrices\footnote{There exists hyperbolicity results for random or deterministic product of matrices like
\cite{bgmv}. However, they are based on genericity arguments which can not be used in order to understand the role of each of
the elements $n$, $\nu$ and $\Pi$ in the creation of hyperbolicity.}. In the following, we solve the hyperbolicity problem
in the three degrees of freedom case.

\section{The transversality-torsion phenomenon}
\label{transtorpheno}

In this section, we deal with three degrees of freedom Hamiltonian systems. In the following, we denote by
${\cal M}_{n,p} (\rR )$ the set of $n\times p$ matrices with real coefficients and
for all matrices $M\in {\cal M}_{n\times n} (\rR )$, we denote by $\mid M \mid$ its determinant.

\subsection{Transversality constraints}
\label{ctrans}

The matrix $\Pi\in {\cal M}_{4,4}(\rR)$ has the following form in the symplectic base $(e_{\phi},e_s, e_{\rho} ,e_u )$:
$$
\Pi =
\left (
\begin{array}{cc}
A & B\\
C & D
\end{array}
\right )
,$$
where $A,B,C,D\in {\cal M}_{2,2}(\rR)$.\\

For all differentiable manifold $\cal M$, we denote by $T_x {\cal M}$ the tangent space to $\cal M$
at point $x\in {\cal M}$.

\begin{defi}
We say that the homoclinic matrix is transverse if and only if it satisfies the following transversality condition
$\Pi (T_{p^-} W^- (T) ) + T_{p^+} W^+ (T) =T_{p^+} {\cal S}$.
\end{defi}

Of course, if the intersection of the stable and unstable manifold, $W^+ (T)$ and $W^- (T)$, of a torus $T$ is transverse
along an homoclinic orbit $\gamma$, then the homoclinic matrix satisfies the transversality conditions by definition.

\begin{lem}
\label{mtrans}
The matrix $\Pi$ is transverse if and only if
$\Delta =\di
\left |
\begin{array}{llll}
c_{1,1} & d_{1,2} \\
c_{2,1} & d_{2,2} \\
\end{array}
\right |
\neq 0$.
\end{lem}

\begin{proof}
Let $v=(v_{\phi} ,0,0,v_u )$ be a vector in $T_{p^-} W^- (T)$. We have
\begin{equation}
\Pi v=(a_{11} v_{\phi} +b_{12} v_u , a_{21} v_{\phi} +b_{22} v_u ,
c_{11} v_{\phi} +d_{12} v_u , c_{21} v_{\phi} +d_{22} v_u ) .
\end{equation}
We begins with the global condition of transversality, namely that $v'=\Pi v=(v'_{\phi} ,
v'_s ,v'_{\rho} ,v'_u )$ is such that $v'_{\rho} =0$ and $v'_u =0$ if and only if
$v_{\phi } =0$ and $v_u =0$. This condition implies
$$\left | \begin{array}{llll}
c_{1,1} & d_{1,2} \\
c_{2,1} & d_{2,2} \\
\end{array}
\right |
\neq 0.$$
\end{proof}

In the following, we need the following strengthening of the transversality condition:

\begin{defi}
The matrix $\Pi$ is strongly transverse if $\Delta\not=0$ and $d_{2,2} \not= 0$.
\end{defi}

The condition $d_{2,2} \not= 0$ does not come from the transversality assumption. We can understand the geometrical
nature of this condition as follow:\\

The unstable (resp. stable) manifold $W^u (T)$ (resp. $W^s (T)$) is foliated by $1$ dimensional manifolds (see \cite{wi},p.138) denoted by
$W^u_p (T)$ (resp. $W^s_p (T)$), $p\in T$ (the Fenichel fibers), and
\begin{equation}
W^u (T)=\bigcup_{p\in T} W^u_p (T) \ \ (\mbox{\rm resp.}\ W^s (T)=\bigcup_{p\in T} W^s_p (T)\, ).
\end{equation}
In the normal form coordinates system, we have for all $p=(\phi_p ,0,0,0) \in T$,
\begin{eqnarray}
W^u_{(\phi_p ,0,0,0)} (T)=\{ (\phi ,s,\rho ,u)\in \T \times \R \times \R\times \R \mid \, \phi=\phi_p ,\
s=0,\ \rho =0 \} ,\\
W^s_{(\phi_p ,0,0,0)} (T)=\{ (\phi ,s,\rho ,u)\in \T \times \R \times \R\times \R \mid \, \phi=\phi_p ,\
u=0,\ \rho =0 \} .
\end{eqnarray}
The condition $d_{2,2}\not= 0$ is then equivalent to the following geometrical condition on the foliation of the stable
and unstable manifolds in the linear model.

\begin{lem}
Let us consider the linear model. The condition $d_{2,2}\not= 0$ is equivalent to the transversality of the intersection
between the unstable leave at $(\phi_- ,0,0,0)$ denoted by $W^u_{(\phi_- ,0,0,0)} (T)$ with the invariant manifold defined
by $\{ (\phi, s,\rho ,u) \in \T \times \R \times \R \times \R ;\, u=0\}$ at point $(\phi_+ ,s_+ ,0,0)$.
\end{lem}

\subsection{The transversality-torsion phenomenon}

The main technical result of this paper is the following:

\begin{thm}[Transversality-torsion phenomenon]
\label{main}
Let $H$ be a 3 degrees of freedom Hamiltonian system  possessing a 2 dimensional partially hyperbolic tori with a transversal
homoclinic connection. We keep notations from section \ref{mise}. We assume that:\\

i) The homoclinic matrix $\Pi$ is transverse;

ii) The torus is with torsion;\\

Then, for $n$ sufficiently large, the matrix $\Pi Df_l^n$ is hyperbolic.

Moreover, if the matrix $\Pi$ is strongly transverse, i.e. $d_{22}\not= 0$, all its eigenvalues are reals and given
asymptotically by
$$x_1 \sim -n\nu d_{22}^{-1} \Delta ,\ x_2 \sim d_{22} \lambda^{-n}
,\ x_3 =x_1^{-1} ,\ x_4 =x_2^{-1} .$$
\end{thm}

\begin{proof}
Let us assume that the matrix $\Pi Df_l^n$ possesses a complex eigenvalue $\beta$. As $\Pi Df_l^n$ is symplectic,
we know that the three remanning eigenvalues are $\bar{\beta}$, $1/\beta$ and $1/\bar{\beta}$ (see
\cite{kh}, prop. 5.5.6, p. 220). The characteristic polynomial is then given by
$$P_n (x)= x^4 +A(n) x^3 +B(n) x^2 +A(n) x +1,$$
where $A(n)= -(S +\bar{S} )$, $B(n)= 2+\mid S\mid^2$ with $S =\beta +\di {1\over \beta}$.

Moreover, we have
$$
\left .
\begin{array}{lll}
A(n) & = & -d_{22} \lambda^{-n} -\lambda^n a_{22} -n\nu c_{11} -a_{11} -d_{11} , \\
B(n) & = & \lambda^n \di\left [ \mid A\mid +a_{22} d_{11}-c_{12} b_{21} +n\nu (a_{22} c_{11} -c_{12} a_{21} ) \di\right ] \\
 & & + \lambda^{-n} \left [ \mid D\mid +a_{11} d_{22} -c_{21} b_{12} +n\nu \Delta \right ] \\
  & & + ( a_{11}d_{11} +a_{22}d_{22} -c_{22} b_{22} -c_{11} b_{11} ) .
\end{array}
\right .
$$
We must consider two cases: $d_{22} \not= 0$ and $d_{22} =0$.\\

- If $d_{22} \not=0$, i.e. we have for $n$ sufficiently large $A(n) \sim -d_{22} \lambda^{-n}$. In the same way, as
$\Delta \not=0$ and $\nu\not= 0$, we obtain $B(n)\sim n\nu \Delta \lambda^{-n}$. We deduce that
$\mbox{\rm Re} S \sim d_{22} \lambda^{-n}$ and $\mid S\mid^2 \sim d_{22}^2 \lambda^{-2n}$. We also have $\mid S\mid ^2
\sim n\nu \Delta \lambda^{-n}$ using the inequality on $B(n)$. We obtain a contradiction. As a consequence,
all the eigenvalues are reals.

We then have eigenvalues $x_1$, $x_2$ and $1/x_1$, $1/x_2$, $x_1 \in \rR$ and $x_2 \in \rR$.
We denote by $S_1 =x_1 +1/x_1$ and $S_2 =x_2 + 1/x_2$.
We have $A(n)=-(S_1 +S_2 )$ and $B(n)=2+S_1 S_2$, so $S_1 (A(n) +S_1 )=-S_1 S_2$. As
$A(n)\sim -d_{22} \lambda^{-n}$ and $B(n) \sim n\nu \Delta \lambda^{-n}$, we conclude that
$S_1 \sim -n d_{22}^{-1} \Delta$, so $x_1 \sim -nd_{22}^{-1} \Delta$. Using
$A(n)$, we obtain $S_2 \sim d_{22} \lambda^{-n}$, so $x_2 \sim d_{22} \lambda^{-n}$, which concludes the proof.\\

- If $d_{22} =0$, we have $A(n) =O (n)$. As $B(n)\sim n\nu \Delta \lambda^{-n}$, this implies that $\mbox{\rm
Im} (S)\not= 0$. If we denote $\beta =\beta_1 +i\beta_2$, $\beta_1 ,\beta_2 \in \R$, we have $\mbox{\rm Im} (S)=\beta_2
(1-\mid \beta \mid^{-1} )$. As $\mbox{\rm Im} (S)\not = 0$, we deduce that $\beta_2 \not= 0$ and $1-\mid \beta\mid^{-1} \not =0$,
i.e. $\beta_2 \not =0$ and $\mid \beta \mid \not= 1$. The eigenvalues are then hyperbolic.\\

This concludes the proof.
\end{proof}

The behaviour of the eigenvalues can be also given when $d_{22} =0$, but depends on several assumptions on the form
of the homoclinic matrix which have not a direct geometrical meaning.\\

In some cases of interest, we can obtain a stronger result. For example, using the homoclinic matrix introduced in \cite{ma}
and generalized in \cite{cr2} in relation with the Arnold model \cite{ar}, we obtain:

\begin{thm}
Let $H$ be a three degrees of freedom Hamiltonian systems possessing a 2 dimensional partially hyperbolic tori with
a transversal homoclinic connection. We keep the notations of $\S$.\ref{mise}. We assume that the homoclinic matrix has the form
\begin{equation}
\Pi =\left (
\begin{array}{cccc}
1 & 0 & 0 & 0 \\
0 & 1 & 0 & 0 \\
\delta & 0 & 1 & 0 \\
0 & 0 & 0 & 1
\end{array}
\right ) ,
\end{equation}
and $\delta$ is a parameter.

Then, the matrix $\Pi Df_l^n$ is hyperbolic for $n$ sufficiently large if and only if the matrix $\Pi$ is transverse,
i.e. $\delta\not= 0$ and the torus is with torsion, i.e. $\nu\not=0$.
\end{thm}

\begin{proof}
The characteristic polynomial of  $\Pi Df_l^n$ is given by
$$P(x)=(x^2 -x(\delta n\nu +2) +1 )(x^2 -xa(n)+1) ,$$
where $a(n)=\lambda^{2n} +\lambda^{-n}$. The matrix is hyperbolic if $\delta\not= 0$ and $\nu\not= 0$. Indeed, in this case,
the matrix $\Pi$ satisfies the transversality assumption. Moreover, if
$\delta=0$ and $\nu\not=0$ (or $\delta\not= 0$ and $\nu=0$), we obtain two eigenvalues equal to $\pm 1$, destroying the
hyperbolicity. This concludes the proof.
\end{proof}

\end{document}